\documentclass[12pt]{article}
\usepackage{amssymb,amsmath}

\title{Non-commutative compactifications and elliptic curves}

\author {Yan Soibelman, Vadim Vologodsky}

\begin{document}
\maketitle

\newtheorem{thm}{Theorem}
\newtheorem{lmm}{Lemma}
\newtheorem{dfn}{Definition}
\newtheorem{rmk}{Remark}
\newtheorem{prp}{Proposition}
\newtheorem{conj}{Conjecture}
\newtheorem{exa}{Example}
\newtheorem{cor}{Corollary}
\newtheorem{que}{Question}
\newtheorem{ack}{Acknowledgements}
\newcommand{\K}{{\bf k}}
\newcommand{\C}{{\bf C}}
\newcommand{\R}{{\bf R}}
\newcommand{\N}{{\bf N}}
\newcommand{\Z}{{\bf Z}}
\newcommand{\Q}{{\bf Q}}
\newcommand{\G}{\Gamma}
\newcommand{\A}{A_{\infty}}

\newcommand{\ihom}{\underline{\Hom}}
\newcommand{\ra}{\longrightarrow}
\newcommand{\epi}{\twoheadrightarrow}
\newcommand{\mono}{\hookrightarrow}

\newcommand{\epp}{\varepsilon}

\centerline{\it To Yuri Ivanovich Manin on his 65th birthday}

\vspace{5mm}

{\bf Abstract.} We introduce non-commutative degenerations of elliptic curves
over local fields. Corresponding objects are close relatives of non-commutative
tori of Connes and Rieffel.

\vspace{10mm}

\section{Introduction}

Let $X$ be an analytic space over a local field ($\C$, $\R$ or $\Q_p$), and $G$
be a group acting freely on $X$.  If the action is not discontinuous the
quotient $X/G$ is not an   analytic  space. Nevertheless one can define the ``category of coherent
sheaves on $X/G$'' (notation $Coh(X/G)$)
 as the category  of $G$-equivariant coherent sheaves on $X$.
A typical example studied in this paper is $X=\C^{\ast}$ and  $G=q^{\Z}$ acting
by dilations $z\mapsto q^nz$. If $|q|\ne 1$ the quotient $\C^{\ast}/q^{\Z}$
is an elliptic curve. If $|q|=1$ the quotient is not even a Hausdorff space.
In either case the category $Coh(X/G)$ is equivalent to the category
of certain modules over the cross-product of the algebra
${\cal O}(\C^{\ast})$ of analytic functions  and the group
algebra of $q^{\Z}$ (see Section 3).
According to the philosophy of non-commutative compactifications
proposed  in [So1]  the compactification of the moduli space of elliptic curves
should  also contain objects corresponding to the case $|q|=1$.
They can be viewed as non-commutative degenerations
of ordinary elliptic curves and called non-commutative (or  quantum) elliptic curves.
They also can be called analytic quantum tori because
in non-commutative geometry $C^{\infty}$-versions of these objects are called
quantum or non-commutative tori. The latter were introduced by Alain Connes and
studied by Connes and Marc Rieffel in the early 80's (see for ex. [Co]).

Elliptic curves (commutative or not) can be viewed as
$q$-difference analogs of analytic vector bundles with  flat connections.
In this sense they are a special case of holonomic $q$-D-modules (cf. [Sab]).
There is a vast literature on $q$-difference equations. It is mainly devoted
to the case $|q|\ne 1$.

From the point of view of non-commutative geometry our non-commuta-\\
tive elliptic curves are ``bad'' quotients. On the other hand, a complex elliptic curve is
a $C^{\infty}$-manifold (2-dimensional torus) equipped with a complex structure.
Hence one has another degeneration: complex structure degenerates into a foliation.
The corresponding categories of vector bundles which carry a flat connection were
studied in [So1]. It would be interesting to make more precise the relationship
between our non-commutative elliptic curves and quantum tori from [So1].
We remark that the picture of [So1] does not have a p-adic analog, while 
non-commutative elliptic curves exist over non-archimedean local fields as well.

About the contents of the paper. In Section 2 we prove a theorem characterizing algebraically
the category of coherent sheaves on certain quotients of Stein varieties.
In Section 3 we define degenerate elliptic curves (there $|q|=1$).
This includes both analytic and formal versions. In particular, in Section 3.2 we revisit
the results of [BG] concerning the formal case.
We compute their Picard groups and obtain some results about irreducible
objects in the category of coherent sheaves. 
We also discuss the analog
of the uniformization of an elliptic curve by a complex line. Using this picture
we define the moduli space of non-commutative elliptic curves and prove
the analog of Rieffel's theorem about $SL_2(\Z)$-action on the moduli
space of quantum tori. Our approach should work in the case of abelian
varieties as well. In this case we expect that the analog of the Rieffel-Schwarz theorem
holds (see [RS]).
The reader will notice that in Section 3 we do not distinguish between ordinary and
degenerate elliptic curves calling all of them non-commutative.
In the last Section 4 we sketch the case of non-arhimedean quantum tori.
There is an evidence that using the ideas of [KoSo1] one can apply them
to quantization of hyper-K\"ahler manifolds.

We do not discuss in the paper the general philosophy of non-commutative
compactifications as well as various connections
and ramifications of this concept. Interested reader can consult with [So1].
There are many speculative relations of non-commutative compactifications
with some structures in number theory, dynamical systems and physics
(for ex. with [BDar], [CDS], [De1-2], [Mo1-2]). We hope to address some of them
in the future.

{\it Acknowledgements.} We thank to Ofer Gabber, Victor Ginzburg,
 Maxim Kontsevich, Yuri Manin,
for interesting discussions. Both authors are thankful to IHES for hospitality
and excellent reasearch atmosphere.
Y.S. gratefully acknowledges  financial support him as a Fellow of the Clay Mathematics Institute.

\section{Coherent sheaves on elliptic curves}

Let $K$ be a either a field of complex numbers
or non-archimedean field, $q\in K^{\ast}, |q|<1$.
Recall that the smooth analytic space $E_q(K)=G_m(K)/q^{\Z}=K^{\ast}/q^{\Z}$
is an elliptic curve (in particular it is an algebraic
variety). We denote by $Coh(E_q)$ the category of
algebraic coherent sheaves on $E_q$. Since
$G_m(K)=K^{\ast}$ is an analytic space
 we have a sheaf of analytic functions ${\cal O}_{K^{\ast}}$,
and the category of analytic coherent sheaves $Coh(K^{\ast})$.
The natural projection $\pi:K^{\ast}\to E_q$ gives rise to a faithful exact
functor (pull-back) $\pi^{\ast}:Coh(E_q)\to Coh(K^{\ast})$.
The image ${\cal C}_q=\pi^{\ast}(Coh(E_q))$ is the category
of $q^{\Z}$-equivariant analytic coherent sheaves on $K^{\ast}$.

The following proposition is well-known (see [Ki] for non-archimedean version).

\begin{prp} If $F\in Coh(K^{\ast})$ then $H^i(K^{\ast},F)=0$
for $i>0$, and each $F_x, x\in K^{\ast} $ is generated by $H^0(K^{\ast},F)$.

\end{prp}

In particular the functor $\phi:V\mapsto H^0(K^{\ast},\pi^{\ast}(V))$
is faithful and exact from the category $Coh(E_q)$ to the category
of modules over the algebra $A_q^{an}$ generated by the analytic functions
$f\in {\cal O}(K^{\ast})$ and the invertible generator $\xi$ such that
$\xi f(z)=f(qz) \xi$ for all $z\in K^{\ast}$. Indeed the space of global sections of an equivariant
sheaf on $K^{\ast}$ gives rise to an $A_q^{an}$-module.
Notice that $A_q^{an}$ contains two commutative subalgebras:
the algebra ${\cal O}(K^{\ast})$ of analytic functions and the
algebra $K[q^{\Z}]:=K[\xi,\xi^{-1}]$.

The aim of this section is to prove the following result.

\begin{thm} The functor $\phi$ provides an equivalence of
$Coh(E_q)$ with the category ${\cal B}_q$ of such $A_q^{an}$-modules $M$
that $M$ is finitely presentable over ${\cal O}(K^{\ast})$,
i.e. there exists $m,n>0$
and an exact sequence of ${\cal O}(K^{\ast})$-modules

$${\cal O}^n(K^{\ast})\to {\cal O}^m(K^{\ast})\to M\to 0.$$

\end{thm}

This theorem can be generalized to abelian varieties, and furthermore
to certain class of Stein spaces. We would like to state and
prove such a generalization.

Let $X$ be a Stein space over $K$, and $G$ be a group freely acting
on $X$. We assume that $X/G$ is a projective algebraic
variety and the pull-back $\pi^{\ast}({\cal O}_{X/G}(1))$ is
a trivial analytic line bundle on $X$ (here $\pi:X\to X/G$
is the natural projection).
The functor $\pi^{\ast}$ gives rise to an equivalence
of the category $Coh(X/G)$ of algebraic coherent sheaves on $X/G$
and the category $Coh_G(X)$ of $G$-equivariant analytic
coherent sheaves on $X$.
The functor of global sections
$\Gamma_X:F\to \Gamma(X,F)$ is an exact functor from $Coh_G(X)$
to the category of modules over the algebra
${\cal O}(X,G):={\cal O}(X)\rtimes G$
(cross-product of $G$ and the algebra
of analytic functions ${\cal O}(X)$).
We would like to describe the image of $Coh(X/G)$ under the functor
$\phi=\Gamma_X\circ \pi^{\ast}$.
In order to do that we introduce the category ${\cal B}_G(X)$
of ${\cal O}(X,G)$-modules $M$ which are
finitely presentable over ${\cal O}(X)$. This means that
there exist $m,n>0$ and an exact sequence of ${\cal O}(X)$-modules

$${\cal O}^n(X)\to {\cal O}^m(X)\to M\to 0.$$

\begin{thm} The functor $\phi$ gives rise to an equivalence
of categories
$$Coh(X/G)\simeq {\cal B}_G(X)$$

\end{thm}

Clearly the previous theorem follows from this one for $X={K}^{\ast}$,
$G={\Z}$ acting by $n(x)=q^nx, n\in {\Z}, x\in {K}^{\ast}.$

{\it Proof.} The proof will consists of several steps.

{\it Step 1.} Let us check that if $F\in Ob(Coh(X/G))$ then
$\phi(F)\in Ob({\cal B}_G(X))$. It follows from the Serre's theorem
that there exist
vector bundles $W_1,W_2$ and a morphism $W_1\to W_2\to F\to 0$.
Moreover, we  can choose $W_k=\oplus_i {\cal O}(m_i), k=1,2$.

Pulling back $W_i, i=1,2$ to $X$
we get a $G$-equivariant analytic bundles $\pi^{\ast}(W_i), i=1,2$
and $G$-equivariant analytic coherent sheaf $\pi^{\ast}(F)$.
The functor $\phi$ is exact as the composition of two exact
functors ($\Gamma_X$ is exact because $X$ is Stein).
Thus we have an exact sequence of ${\cal O}(X)$-modules

$$\phi(W_1)\to \phi(W_2)\to \phi(F)\to 0.$$

By our assumption
the pull-backs of the sheaves ${\cal O}_{X/G}(l), l\in \Z$
are trivial analytic vector bundles on $X$.
Hence $\phi(W_i), i=1,2$ are free  ${\cal O}(X)$-modules
of finite ranks. It follows that
$\phi(F)$ is an object of the category ${\cal B}_G(X).$

{\it Step 2.} Observe that if $F\ne 0$ then $\phi(F)\ne 0$.
Indeed,  $\pi^{\ast}(F)\ne 0$ because $\pi^{\ast}$ is a faithful
functor. Since $X$ is Stein the sheaf $\pi^{\ast}(F)$ is generated
by its global sections which is $\phi(F)$. Hence if $\phi(F)=0$
then $F=0$. Similarly one proves that if $f:V\to W$ is a
non-trivial morphism
of coherent sheaves on $X/G$ then $\phi(f)\ne 0$.
Indeed, consider the subsheaf $Im(f)$ of $ W$. Then
$\pi^{\ast}(Im(f))$ is generated by $\Gamma(X,\pi^{\ast}(Im(f)))\simeq
Im(\phi(f)).$

Thus we have proved that $\phi:Coh(X/G)\to {\cal B}_G(X)$ is an
exact faithful functor.

{\it Step 3.} We would like to prove that any object
of ${\cal B}_G(X)$ is isomorphic to some $\phi(F), F\in Coh(X/G).$

Let $M$ be an object of the category  ${\cal B}_G(X)$.
Then for some $m,n>0$
we have an exact sequence
of ${\cal O}(X)$-modules

$${\cal O}^n(X)\to {\cal O}^m(X)\to M\to 0.$$

We consider the localization functor
$M\mapsto \widetilde{M}:={\cal O}_X \otimes M$, where ${\cal O}_X$
is the sheaf of analytic functions on $X$.
The localization functor is not faithful on finitely generated
modules. It is easy to construct an example of a module $V\ne 0$
such that $\widetilde{V}=0$. Nevertheless we have the following result.

\begin{lmm} The localization functor gives rise to
an equivalence of the category  of finitely representable  ${\cal O}(X)$-modules
and the full subcategory of coherent sheaves whose objects are
sheaves isomorphic to $coker({\cal O}^n_X \stackrel{f}{\to}{\cal O}^m_X)$
for some $n$, $m$ and $f$.

The inverse functor $\Gamma_X$ associates to a sheaf $F$ the
${\cal O}(X)$-module
of its global sections.
\end{lmm}
The lemma follows from the exactness of the localization functor
and the functor $\Gamma_X$ and from fact that
$\widetilde {{\cal O}^n(X)}= {\cal O}^n_X$, $\, \Gamma_X({\cal O}^n_X)=
{\cal O}^n(X)$.

Returning to the proof of the theorem
we see that
   $$ \Gamma_X( \widetilde{F})= F$$
for any $F\in Coh(X/G)$. Summarizing, we have constructed an exact fully faithful
functor $\phi$ between the categories in question.
This completes the proof. $\blacksquare$

We do not know whether Theorem 2 holds without
 the assumption of triviality of the pull-back of   ${\cal O}_{X/G}(1)$.

\section {Elliptic curves for $|q|=1$}

\subsection {Category ${\cal B}_q$}

Let $K$ be either the field $\C$ or a non-archimedean local field.
Let $X$ be  Stein space over $K$, ${\cal O}(X)$ algebra
of analytic functions on $X$,  $G$ a group acting freely on $X$.

\begin{dfn} The category of coherent sheaves on the non-commutative
space $X/G$ is defined as a category of modules over the cross-product algebra
${\cal O}(X)\rtimes G$ which are finitely presentable over ${\cal O}(X)$.

\end{dfn}

From the point of view of this definition one should not distinguish between
``good'' and ``bad'' quotients. They all are non-commutative spaces.

Let us assume that $K={\C}$.
We will use the above definition in two cases: $X=\C^{\ast}, G=q^{\Z}$ and
$X=\C, G=\Z\oplus \Z$ with the natural actions of $G$
(by dilations in the first case and by shifts by the lattice
in the second one). In this section we will study the case $|q|=1$ or when
the action of the lattice on $\C$ is not discontinuous.

By the very definition the category ${\cal B}_1$ is the category 
of coherent sheaves on ${\C}^{\ast}$, whose space of global sections is 
finitely presentable over ${\cal O}({\C}^{\ast})$, endowed with an automorphism.

It is likely that the category ${\cal B}_q$ is
also equivalent to the latter,  provided that $q$ is a root of unity. We do not prove this fact
here\footnote{Reader may wish to look at the beginning of Section 3.2 where a similar result is proven.} and will not use it.
  
From now on we assume that  $|q|=1$ and $q$ is not a root of $1$.
We start with the following lemma.

\begin{lmm}
 Let $V\in Ob({\cal B}_q)$. Then the corresponding
coherent sheaf $\widetilde{V}$  is  free of finite rank.

\end{lmm}

{\it Proof.} The sheaf  $\widetilde{V}$ is isomorphic to
 $coker({\cal O}^n_{{\C}^{\ast}} \stackrel{f}{\to}{\cal O}^m_{{\C}^{\ast}})$
for some $n$, $m$ and $f$. Let $k$ be the maximal rank of
$f_x:{\C}^n \to  {\C}^m$ for $x\in {\C}^{\ast}$
and $S$ be the subset of ${\C}^{\ast}$ where the
rank is strictly smaller the $k$. Since $S$ is contained in the set of zeros
of a holomorphic function on  ${\C}^{\ast}$ (a minor of the matrix $f$), it is
discrete. On the other hand, $S$ must be  $\Z$-invariant. Hence $S$ is empty
i.e.   $\widetilde{V}$ is a vector bundle. It remains to observe that
any holomorphic vector bundle on ${\C}^{\ast}$ is trivial.

\begin{prp} The category ${\cal B}_q$ is abelian for all $q\ne 0$.

\end{prp}

{\it Proof.}
 For any object $V$ of ${\cal B}_q$ we have the corresponding
coherent sheaf $\widetilde{V}$ on ${\C}^{\ast}$. If $f:V\to W$
is a morphism in ${\cal B}_q$ then $\widetilde{Ker(f)}$
and $\widetilde{Coker(f)}$ are coherent
sheaves. According to the previous lemma they are both free
sheaves of finite rank. Hence the spaces of sections
$Ker(f)$ and $Coker(f)$ are $A_q^{an}$-modules which are
of the type ${\cal O}^p({\C}^{\ast})$ as ${\cal O}({\C}^{\ast})$-modules.
Hence they belong to ${\cal B}_q$. $\blacksquare$

Let $K({\C}^{\ast})_{mer}$ be the field of meromorphic functions
on ${\C}^{\ast}$. Then we have a functor $\mu:V\mapsto
V\otimes_{{\cal O}({\C}^{\ast})} K({\C}^{\ast})_{mer}$
from the category ${\cal B}_q$ to the category
of finite-dimensional $K({\C}^{\ast})_{mer}$-vector spaces.

\begin{prp} The functor $\mu$ is fully faithful.

\end{prp}

{\it Proof.}
Every object of ${\cal B}_q$ gives rise to a coherent sheaf
on  ${\C}^{\ast}$. Hence it suffices to prove the corresponding statement
for coherent sheaves. It is easy to see that if $V\ne 0$ then
$\mu(V)\ne 0$ (otherwise $V$ is a torsion analytic sheaf which
is $q^{\Z}$ equivariant. This is impossible because orbits
of the $q^{\Z}$-action are dense in circles $|z|$=const).
Applying this observation to the kernel of a morphism we see that
$\mu$ is a monomorphism on morphisms.
 $\blacksquare$

 Notice that the functor $\xi$ is not fully faithful  in the case $|q|\ne 1$:
it kills all torsion sheaves.

\begin{rmk} The above results hold for $K=\Q_p$ as well.

\end{rmk}

\begin{prp} Suppose that
$q$ does not satisfy the Liouville property, i.e. there exists $L>0$ such that
$|q^n-1|=O(n^{-L})$. Then the Picard group  (i.e. the
 of isomorphism classes of
 objects of  ${\cal B}_q$
which are rank one modules over  ${\cal O}(\C ^{\ast})$) is
isomorphic to an extension of $\Z$ by ${\C}^{\ast}/q^{\Z}$.

\end{prp}

{\it Proof}.  Let  $A={\cal O}(\C ^{\ast})$. Then $A$ is a $\Z$-module:
the generator $1\in \Z$ sends $f(z)$ to $f(qz)$.
Clearly the Picard group is isomophic to $H^1(\Z, A^{\ast})$.
The proposition follows easily from the exact sequence
$$0 \to \Z \to A \stackrel{exp}{\longrightarrow} A^{\ast}
\stackrel{v}{\longrightarrow} \Z\to 0$$
(here $v(f)=\int d\,log\,f=0$ where the integral
is taken over any small circle around $z=0$) and from the following
 lemma

\begin{lmm} Suppose that
$q$ does not satisfy the Liouville property. Then  $H^1(\Z, A)=\C$
\end{lmm}

{\it Proof}. Let $g(z)=\sum^{+\infty}_{-\infty}g_n z^n$.
If $g_0=0$ then the equation $g(z)=a(qz)-a(z)$ has a solution
$a(z)=\sum^{+\infty}_{-\infty}a_n z^n$,
where $a_0=0, a_n=g_n/(q^n-1), n\ne 0$. If $q$ does not satisfy the Liouville property
the series converges. $\blacksquare$

\begin{rmk} It follows from the proof above that for any $q$, which is not a root of $1$
 there is a canonical
surjective homomorphism:
     $$deg:\,  Pic({\cal B}_q) \to \Z $$

The number $deg(V)$ is called the degree of the object $V$.

\end{rmk}

It is an interesting problem to describe simple objects
of the category ${\cal B}_q, |q|=1$.
In the case $|q|\ne 1$ it is equivalent to finding all
irreducible sheaves on the elliptic curve $E_q$.
They correspond to points of $E_q$ (torsion sheaves).
A vector bundle is not a simple object in the category
of coherent sheaves.
The situation changes in the case $|q|=1$.
It follows from Lemma 3 that if $|q|=1$ and $q$ is not a root of $1$,
 there are no ``torsion sheaves" at all.
On the other hand,  objects of rank $1$
over  ${\cal O}(\C ^{\ast})$ are simple.

We are going to construct some simple objects of higher ranks.
For general $q$ and $s$ with $s^n=q$ we have the obvious $n$-sheets coverings
$f_n: \C ^{\ast}/s^{\Z}\to \C ^{\ast}/q^{\Z}$ and
$g_n: \C ^{\ast}/q^{\Z}\to \C ^{\ast}/s^{\Z}$.
The coverings give rise to the functors
$$f_{n \ast}, g_n^{\ast}: {\cal B}_s \to {\cal B}_q \;
\;, g_{n \ast}, f_n^{\ast}: {\cal B}_q \to {\cal B}_s $$

Mimicking the usual construction one can define the above functors
for {\it any} $q\in \C^{\ast}$. The functor $f_{n\ast}$ corresponds to the following
homomorphism of algebras
  $$A^{an}_q \to A^{an}_s, \; f(z) \mapsto f(z^n) , \;  \xi \to \xi, $$

while the functor $g_n^{\ast}$ corresponds to the homomorphism

  $$A^{an}_q \to A^{an}_s, \; f(z) \mapsto f(z) , \;  \xi \to \xi ^n. $$

The other two functors are defined as the tensor product
 $A^{an}_q \otimes _ {A^{an}_s } V$, where in the case  of $f_{n\ast}$ the
$A^{an}_s$-module structure on    $A^{an}_q $ is given by
 the first above homomorphism of algebras while in the case of
 $g^{\ast}_n$ it is given by the second one.
One can prove that pull-back and push-forward functors satisfy the
usual properties.

Remind (see the proof of Proposition 4) that there is a canonical morphism
from ${\cal O}(\C ^{\ast})$ to the Picard group of ${\cal B}_s$.
Given integers $k$ and $n$ with $n>0$ we denote by
$P_{k,1}$ the rank one module corresponding to the function
$z^k$ and by  $P_{k,n}$ the push-forward $f_{n\ast}P_{k,1}$.
The number $k$ coincides with the degree $deg(P_{k,n})$.

\begin{lmm}
Assume that $|q|=1$ and $q$ is not a root of unity.
Then, for any integer $k$ and $n>0$ such that $(k;n)=1$, the module  $P_{k,n}$ is simple.
\end{lmm}

{\it Proof.} Assume the contrary, and  choose a proper submodule $W\subset P_{k,n}$.
Consider the pull back $f^{\ast}_n (W)\subset f^{\ast}_n  (P_{k,n})$. It is easy to show that
$ f^{\ast}_n  (P_{k,n})$ is isomorphic to the direct sum of rank $1$ objects of degree $k$,
namely

$$  f^{\ast}_n  (P_{k,n})= P_{k,1}\otimes (\oplus _{\mu,\, \mu^n=1}A^{an}_s/(\xi -\mu)),$$
where as usual we denote by $(\xi -\mu)$ the right ideal generated by the expression in the brackets.

Since all objects $P_{k,1}\otimes A^{an}_s/(\xi -\mu)$ are simple and pairwise
 non-isomorphic it follows that
$f^{\ast}_n (W)$ is isomorphic to the  direct sum of  rank $1$ objects of the form
$P_{k,1}\otimes A^{an}_s/(\xi -\mu)$. In particular,
$deg (det (f^{\ast}_n (W)))= kl$, where
$0 <l=rk (W) < n$. On the other hand, $deg (det (f^{\ast}_n W))=
deg(f^{\ast}_n(det\,W))= n\, deg(det\,W)$. This contradiction completes the proof.$\blacksquare$

Complete description of simple objects is not known.

\begin{que} How to describe all simple objects of ${\cal B}_q$
in the case $|q|=1$?

\end{que}

At the end of this subsection we will make few comments about
the tensor structure on the category ${\cal B}_q$ (for arbitrary ground field $K$).
Let us fix $q\in K^{\ast}$.
Notice that  the category
of $A^{an}_q$-modules is a $K$-linear symmetric monoidal category.
It follows from the interpretation
of objects of ${\cal B}_q$ as coherent sheaves that the tensor product of
finitely presentable over ${\cal O}(K^{\ast})-$ modules is again a finitely
presentable module. Hence  ${\cal B}_q$ is  a symmetric monoidal category.
It can be considered as a difference analog of the category of analytic vector
bundles on $K^{\ast}$ equipped with a flat connection.

\subsection{Formal analog of the category ${\cal B}_q$}

Let $F:= \C((z))$ and $R:=\C[[z]]$. We denote by  ${\cal B}^{al}_q$ the category of
 finite-dimensional vector spaces over $F$ equipped with a semi-linear
automorphism $\xi : V\to V $, $\xi f(z) v = f(qz) \xi v$ . In  other words, this is a category
of finite-dimensional representations of the algebra  $A^{al}_q$  of non-commutative Laurent
 polynomials $F[\xi,\xi^{-1}]$.

\begin{lmm}
Assume that $q$ is a root of unity. Then the category  ${\cal B}^{al}_q$ is equivalent to the category of finite-dimensional
$F$-vector spaces endowed with an automorphism.
\end{lmm}

{\it Proof.} Assume that $q$ is a primitive root of unity of order $n$ (i.e. $q^n=1$ and $q^k\ne 1$, for $k<n$). 
It is easy to see that the center  of the algebra  
$A^{al}_q$ is $\C((z^n))[\xi^n,\xi^{-n}] \subset F[\xi,\xi^{-1}]$. Localizing we make $A^{al}_q$ into a sheaf of algebras over the
$Spec ( \C((z^n))[\xi^n,\xi^{-n}])= \mathbb{G}_m$.  
One can easily check that $A^{al}_q$ is an Azumaya algebra over $\mathbb{G}_m$. An object of ${\cal B}^{al}_q$ gives rise to 
a sheaf of modules over the Azumaya algebra supported on a closed 
subscheme $S\subset \mathbb{G}_m $ finite over $Spec(F)$ .
The Brauer group of $F$ is trivial. Hence the restriction of 
the Azumaya algebra  $A^{al}_q$ to 
the ind-scheme $\varinjlim S$ (the limit is taken over all  finite closed subschemes $S\subset \mathbb{G}_m $)
splits, i.e. isomorphic to the matrix algebra. It remains to notice that the category of modules over a split
 Azumaya algebra is equivalent to the category of ${\cal O}$-modules on the underlying scheme.$\blacksquare$
  
 In what follows we assume that $q$ is not a root of unity.

Following [BG] we call an object $V$ of ${\cal B}^{al}_q$ integral if there exists a
$\C[\xi, \xi ^{-1}]$-invariant  $R$-lattice $L\subset V$. Integral objects
form an abelian subcategory ${\cal B}^f_q$ of  ${\cal B}^{al}_q$.
The following theorem was proved in [BG].

\begin{thm} The category ${\cal B}^f_q$ is equivalent to the category of
finite-dimensional $\C^{\ast}/q^{\Z}$- graded complex vector  spaces endowed with a
nilpotent operator  which preserves the grading.
\end{thm}

Let us  describe a construction of the equivalence functor. It is based on the following
easy lemma, proof of which is left to the reader.

\begin{lmm}
Let $V$ be an object of  ${\cal B}^{al}_q$, $L\subset V$
be a $\C[\xi]$-invariant $R$-lattice (we do not require $L$ to be
$\xi^{-1}$-invariant), and $\overline P \subset L/zL :=  \overline L$ be a
$\C[{\xi}]$-invariant subspace such that for any eigenvalue
$\alpha$ of $ {\xi}: \overline P \to \overline P $ the number
$q^n \alpha $ is not an eigenvalue of ${\xi}$
on $\overline L$ for all  $n\ne 0$. Then the embedding
$\overline P \hookrightarrow \overline L$ has a $\C[\xi]$-equivariant lifting 
$P\hookrightarrow L$.

\end{lmm}

Given an object $V$ of  ${\cal B}^f_q$  we can find a
$\C[\xi, \xi ^{-1}]$-invariant  $R$-lattice $L\subset V$ such that the eigenvalues of
${\xi}$ on $\overline L$ satisfy the property in  the lemma.\footnote{To construct such a lattice
$L$ we can use the following procedure. First choose any  $\C[\xi, \xi ^{-1}]$-invariant  $R$-lattice $L_0$.
Assume, that it does satisfy the eigenvalues property: numbers $\alpha$ and $q^n \alpha$ are eigenvalues, for 
some $\alpha \ne 0$ and $n>0$. Let   $\bigoplus_{\gamma} \overline L^{\gamma}_0=\overline L_0:= L_0/ zL_0$
be the Jordan decomposition of $\overline L_0$ (i.e. $(\xi - \gamma)^N  \overline L^{\gamma}_0=0$, for large 
$N$). Choose a decomposition $L_0 =\bigoplus_{\gamma} L^{\gamma}_0$ of the $R$-lattice which descends to the latter one modulo
$zR$. Put $L_1: = \bigoplus _ {\gamma \ne q^n \alpha } L^{\gamma}_0 \oplus z^{-1} L^{ q^n \alpha}_0 $. It is easy to see that
$L_1$ is a ${\C}[\xi, \xi ^{-1}]$-invariant  $R$-lattice. If it still does not   satisfy the property in  the lemma we apply the
above procedure again to $L_1$ and so on.} It follows that
$V$ is isomorphic to $  \overline L \otimes _{\C} F$ as a $F[\xi ]$-module.
The action of $\xi$ on $\overline L$  defines a $\C^{\ast}$-grading
($\overline L _a:= ker (\xi -a)^n$ for large $n$). The latter gives rise to a  $\C^{\ast}/q^{\Z}$-
grading.    Finally, we endow $\overline L _a$ with a nilpotent operator $log(a^{-1}\xi)$.
It gives rise to the equivalence functor.

Notice that the argument above shows that there is a functor from the category  ${\cal B}^f_q$ to the
category of finitely generated modules over the algebra of non-commutative  Laurent
 polynomials $\C[z, z^{-1}, \xi,\xi^{-1}] \subset  A^{al}_q$. The latter sends
$V$ to  $\tilde {V} =  \overline L \otimes _{\C} \C[z, z^{-1}] $.

Now we want to obtain a similar description of the whole  category ${\cal B}^{al}_q$.
We claim that  for any object  $V$ of  ${\cal B}^{al}_q$, there exists an integer $n$
such that the pull back $f^*_n (V)$ is a direct sum of  tensor products
of integral objects and one-dimensional ones.  More precisely we have the following result.

\begin{thm}
The category ${\cal B}^{al}_q$ is equivalent to the category of
$\Q$-graded objects of  ${\cal B}^f_q$ (i.e. each object of
 ${\cal B}^{al}_q$ is direct sum of $\Q$ copies of objects of
${\cal B}^f_q$).

\end{thm}

{\it Proof.} Given a rational number $\lambda = \frac{k}{n}$  with  $n>0$, $(n; k)=1$, $k\ne 0$ unless $n=1$,
  we, first, construct a functor
$$G_{k,n}:  {\cal B}^f_q \to {\cal B}^{al}_q$$

For this purpose we represent the category ${\cal B}^f_q$ as a certain ``tensor quotient'' of the category
${\cal B}^f_s$. Namely, consider the full subcategory ${\cal E}$ of ${\cal B}^f_s$, whose objects are
$A^{al}_s$-modules isomorphic to $A^{al}_s/(\xi - \nu) $,
where $\nu \in \C^{\ast}$ and  $\nu ^n=1$. Notice that
${\cal E}$ is a Picard category i.e. a groupoid with a tensor structure such that any object is invertible.
Define the  ``tensor quotient'' ${\cal B}^f_s/{\cal E}$ to be the category whose objects are those
of  ${\cal B}^f_s$ and morphisms between $V$ and $W$ are given by homomorphisms of
$A^{al}_s$-modules
$V\otimes M \to W$ where $M = \bigoplus _{\{\nu|\nu^n=1\}} A^{al}_s/(\xi - \nu) $. This is
an abelian category. For any object $V$ of ${\cal B}^f_s$ and $E$  of   ${\cal E}$ the object
$V\otimes E$ is isomorphic to $V$ . It is easy to see that the functor
   $$g^{\ast}_n: {\cal B}^f_s/{\cal E} \to {\cal B}^f_q$$
is an equivalence.

Given an object $V$ of  ${\cal B}^f_s/{\cal E}$
we define $G_{k,n}(V)= f_{n\,\ast }(V\otimes P_{k,1})$.
(Recall that $ P_{k,1}$ stands for  $A^{al}_s/(\xi - z^k) $). 

We claim that  $G_{k,n}$ extends to morphisms thus gives rise to a functor. Indeed,
$$Hom_{{\cal B}^{al}_q}(f_{n\,\ast }(V\otimes P_{k,1}); f_{n\,\ast }(W\otimes P_{k,1}))=
Hom_{{\cal B}^{al}_s}(f_{ n}^{\ast}f_{n\,\ast }(V\otimes P_{k,1}); W\otimes P_{k,1})$$
Now the claim is an immediate corollary of the following formula:
\begin{equation}\label{formula}
   f_{ n}^{\ast}f_{n\,\ast }(V\otimes P_{k,1})\simeq V\otimes P_{k,1} \otimes M 
\end{equation}
We leave its proof to the reader.

Notice, that the latter argument  also shows that the functor $G_{k,n}$ is  fully faithful.
 
Next, we are going to show that the functor
 $$\oplus_{n,m}G_{k,n}: \bigoplus_{k,n} {\cal B}^f_q \to {\cal B}^{al}_q$$
is an equivalence of categories.

First, let us prove that the functor is fully faithful. Indeed, it remains to show that
$$Hom_{{\cal B}^{al}_q}(G_{k,n}(V), G_{k^{\prime},n^{\prime}}(W))=0 $$
unless $n=n^{\prime}$ and $k=k^{\prime}$. Notice, that the group under consideration is a subgroup of
$$Hom(f_{ nn^{\prime}}^{\ast}f_{n\ast }(P_{k,1}\otimes V);f_{ nn^{\prime}}^{\ast}f_{ n^{\prime}\ast}(P_{k^{\prime},1} \otimes W) )$$
The formula (\ref{formula}) easily implies that the latter is trivial.
 
It remains to show that the image of $\oplus_{n,k}G_{k,n}$ is ${\cal B}^{al}_q$.
\begin{prp}
For any  object $V$ of ${\cal B}^{al}_q$ there exists $n>0$ such that
$f^{\ast}_n (V)$ is a direct sum of $A^{al}_q/(\xi - a z^k)^l$ where $ k, l>0$
are integers and $a\in \C$.
\end{prp}
{\it Proof.} We start with the following lemma.
 \begin{lmm}
All one-sided ideals of  $A^{al}_q$ are principal.
\end{lmm}
{\it Proof of the Lemma}.
 For an element $X\in A^{al}_q$, $X = a_l \xi ^l +...+ a_m \xi ^m $ with $a_l, a_m \ne 0$
we define $s(X) =m-l$. It is easy to see that   $A^{al}_q$ is an Euclidean domain relative to
the function  $s(X)$.

It follows that any finite-dimensional $A^{al}_q$-module is a direct sum of cyclic ones
(i.e. those generated by a single element).
Hence, it is enough to prove the proposition for $V= A^{al}_q/(X)$, where $X\in  A^{al}_q$.
We will do it by induction
on the dimension of $V$.

Without loss of generality we may assume that
 $X= c \xi^l + a_{l-1} \xi ^{l-1} +...+ a_0$ where $a_i$ are Taylor
series in $z$ and $c\in \C^{\ast}$.
 Notice that $L:=Im(R[\xi] \to V)$ is a $\xi$-invariant  $R$-lattice in $V$. If $a_s(0)\ne 0$,
for some $s$, the element $\xi$ induces the operator $\overline {\xi}:L/zL \to L/zL$
which is not nilpotent. Then we can make use of Lemma $6$ to reduce the dimension
of $V$. In general, there exist integers $n$
and $k$ such that  $c(z^{k/n}\xi)^l + a_{l-1}(z^{k/n} \xi )^{l-1} +
...+ a_0= z^{lk/n}(c^{\prime}\xi^l+b_{l-1}\xi^{l-1}+...+b_0)$ with $b_i\in \C[[z]]$ and
$b_s(0) \ne 0$, for some $s$,
and we can apply Lemma $6$ to
$ P_{-k,1}\otimes f^{\ast}_n (V) $.
This completes the proof of the proposition.

Returning to the proof of the theorem we notice that
 that   the modules $A^{al}_q/(\xi - a z^k)^l$ are in the image of $\oplus_{n,k}G_{k,n}$.
The second remark is that the latter image is closed under the push-forward functor
 and under taking subobjects.

Now, given an object $V$ of  ${\cal B}^{al}_q$ we choose $n$
as in the proposition and consider the canonical morphism
  $$ f_{n\ast }f^{\ast}_n (V)\to V $$
It is easy to see that the morphism is surjective.
Combining this with the two observations above we see that
$V$ is in the image of  $\oplus_{n,k}G_{k,n}$. The theorem is proven. $\blacksquare$

\begin{rmk} 1) The argument
given above is  just a variant of  Manin's proof of
the  Classification Theorem for Dieudonne modules (see  [M4]).

2) One can compare our Theorem $4$ with a well-known result from the theory
of $D$-modules which says that,
for any vector bundle $E$ with a connection on the
punctured formal disk, the pull-back of $E$ on some finite covering of the punctured disk is
 a direct sum of  tensor products of D-modules with regular singularities (at the origin)
 and one-dimensional ones.\footnote{We could not find a reference for this result. We are grateful to Roman Bezrukavnikov
who explained it to one of us. He, in turn, refers to the lectures 
by Spencer Bloch at the University of Chicago.}

\end{rmk}

\begin{cor}
Any finite-dimensional  $A^{al}_q$-module
is a direct sum of modules of the type  $A^{al}_q/(\xi^n - a z^k)^l$ where $n>0, k, l>0$
are integers and $a\in {\C}^{\ast}$.
\end{cor}
\begin{cor}
For any two complex numbers $q$ and $q^{\prime}$ (which are not roots of $1$)
the categories  ${\cal B}^{al}_q$
and  ${\cal B}^{al}_{q^{\prime}}$ are equivalent.
\end{cor}

We would like to describe the relation between the category 
${\cal B}^{al}_q$ and its analytic counterpart ${\cal B}_q$. 

First, we construct  a tensor functor
\begin{equation}\label{gb}
{\cal B}^{al}_q \to {\cal B}_q
\end{equation}
In the process of proving Theorem 3 we  observed that there was a functor from the category
${\cal B}^f_q$ to the category of finitely generated modules over  
$\C[z, z^{-1}, \xi,\xi^{-1}]$.
Composing it with the tensor product with $A^{an}_q$
over  the subalgebra $\C[z, z^{-1}, \xi,\xi^{-1}]\subset A^{an}_q$
we obtain the functor ${\cal B}^{f}_q \to {\cal B}_q$.
It immediately follows from Theorem 4
 (and from the construction of the equivalence given in the proof)  that
this functor  extends to (\ref{gb}).

Assume that $|q|\ne 1$. It is easy to see that the image of  (\ref{gb}) lies in the subcategory
$Bun(E) \subset {\cal B}_q$ of vector bundles on $E:=\C^{\ast}/q^{\Z}$.
\begin{prp}
Assume that $|q|\ne 1$. Then functor ${\cal B}^{al}_q \to Bun(E)$
induces a bijection on the set of isomorphism classes of
objects.
\end{prp}
{\it Proof.}
It follows from   the classification  of vector bundles on $E$
 obtained by M. Atiyah  that any  indecomposable vector bundle is of the form 
 $A^{an}_q/(\xi^n - a z^k)^l$ where $n>0, k, l>0$, $(n,k)=1$, $a\in \C^{\ast}$. Therefore, the functor 
induces a bijection on the set of isomorphism classes of indecomposable  objects. 
The proposition immediately follows from this fact. $\blacksquare$.

It was shown in  [BG]  that  the restriction of (\ref{gb}) to 
 the subcategory ${\cal B}^f_q$ gives rise to an equivalence between the latter category and 
the tensor category of semistable vector bundles on $E:=\C^{\ast}/q^{\Z}$ of degree zero.

Of course the functor  ${\cal B}^{al}_q \to Bun(E)$ is not an equivalence, since $Bun(E)$ is not abelian.

We can reformulate Proposition $6$ by saying that there is a bijection between 
the non-abelian  cohomology group $H^1(\Z, GL(n,F)), F=\C((z))$ and the   set of isomorphism 
classes of $n$-dimensional
vector bundles on $E$.

Now, let $G$ be any algebraic group over $\C$. 
Using the Tannakian formalism and Proposition $6$ one can  easily construct a map
form  $H^1(\Z, G(F))$ to the set  of isomorphism classes of principal $G$-bundles on $E$. 
We do not know if it is a bijection.

 According to [BG] this is true if restrict ourselves to the
subset \\ 
$Im(H^1(\Z, G(R))\to H^1(\Z, G(F)))$,
on the one hand, and  semistable $G$-bundles of degree zero, on the other hand. 

If  $|q| = 1$, we can not say much more about  (\ref{gb}), except  that
in  this case ${\cal B}_q$ has more objects then its algebraic counterpart
${\cal B}^{al}_q$. Probably the Corollary 2 is not true for ${\cal B}_q$.

\subsection{Uniformization by the complex line}

We still assume that $K=\C$.
Complex elliptic curves can be also represented as quotients $\C/L$ where
$L$ is a lattice in $\C$. If $L$ is not discrete then the quotient can be considered
as a non-commutative space.
 By definition (cf. [M2]) a {\it quasi-lattice} in $L\subset \C^n$ is a
free abelian subgroup of $\C^n$ of rank $2n$. If $L$ is discrete in the standard
topology then it is called {\it lattice}.

Let $L\subset \C$ be a quasi-lattice.  We define the algebra
$C_L$ as the cross-product of the algebra of analytic functions
${\cal O}({\C})$ and the group algebra of $L$.
It is generated by analytic functions $f(z)$ and  invertible
generators $e_{\lambda}, \lambda \in L$
subject to the relations

$$f_1(z)f_2(z)=f_2(z)f_1(z), e_{\lambda}e_{\mu}=e_{\lambda+\mu}, e_{\lambda}f(z)=
f(z+\lambda)e_{\lambda}.$$

 We define the category
${\cal D}_L$ as the category of $C_L$-modules which are finitely
presentable over ${\cal O}({\C})$. Let $\tau\in \C\setminus \Q,
q=exp(2\pi i \tau)$.  Let us consider a quasi-lattice
$L_{\tau}=\Z\oplus \Z\tau.$ The corresponding non-commutative elliptic curve
will be denoted by $E_{\tau}$.

\begin{prp} The categories ${\cal D}_{L_{\tau}}$ and ${\cal B}_q$
are  equivalent.

\end{prp}

{\it Proof.} The category ${\cal D}_{L_{\tau}}$ can be
identified with the category of $L_{\tau}$-equivariant
analytic coherent sheaves on $\C$. Similarly the category
${\cal B}_q$ can be identified with the category
of analytic $q^{\Z}$-equivariant sheaves on $\C^{\ast}$.
Clearly the category of sheaves are equivalent. $\blacksquare$

\subsection{Moduli space of quantum elliptic curves}

Let $K$ be  a local  field, as in
Section 2.1, $q\in K^{\ast}$.

Recall that the category ${\cal B}_q$ defines a {\it non-commutative (or quantum)}
elliptic curve $E_q$. Of course, if $|q|\ne 1$ it is just the usual elliptic curve over $K$.
The question arises when two quantum elliptic curves should be called equivalent.
It is natural to expect (cf.  [Ga]) that the equivalences are defined by
certain functors between the corresponding
categories ${\cal B}_q$. We do not know such a definition in general.
One of the reasons is lack of the general theory of non-commutative analytic spaces.
Existing results  in non-commutative algebraic geometry (see for ex. [KR]) are not
sufficient for our purposes.
Hovewer, in the case $K=\C$ we can resolve the difficulty
by using the uniformization of elliptic curves
by the complex line.

\begin{dfn}
 We say that that two non-commutative elliptic curves $E_L$ and $E_{L^{\prime}}$
corresponding
to the categories ${\cal D}_L$ and ${\cal D}_{L^{\prime}}$  are equivalent if
there is a continuous isomorphism of
the algebras $C_L$ and $C_{L^{\prime}}$ (the topology on the algebras is induced by the topology
on the algebra of analytic functions ${\cal O}(\C)$).

Moreover we define {\sf groupoid} of non-commutative elliptic curves
putting $Isom(E_L,E_{L^{\prime}})$ to be the set of all continuous
 isomorphisms of
the algebras $C_L$ and $C_{L^{\prime}}$ modulo the group of
inner automorphisms of  $C_L$.
\end{dfn}

\begin{rmk}
The reader familiar with non-commutative
geometry has noticed  that the relationship between the algebras $C_{L}$
and ${A}_q^{an}$ is similar to the one between the
foliation algebra of a torus and the corresponding quantum torus in [Co].
In particular $C_{L}$ gives rise to a generator-free description of non-commutative elliptic curves.
\end{rmk}

Given  a quasi-lattice  $L\subset \C$ we define a group $G_L$ to be the subgroup
of $\C^{\ast}$ consisting of all $\alpha \in \C^{\ast}$ such that $L=\alpha L $.

\begin{thm} a) Non-commutative elliptic curves
 $E_L$ and $E_{L^{\prime}}$ are equivalent if and only if
$L$ and $L^{\prime}$ are similar (i.e. there exists $\alpha\in \C^{\ast}$ such that
$L=\alpha L^{\prime}$).

b) $Isom(E_L,E_{L})$ is isomorphic to the semi-direct product
$(G_L\triangleright \C/L)\triangleright Pic(E_L)$ (clearly $\C/L$ is
a normal subgroup in $G_L$, and  $Pic(E_L)$ is a normal subgroup in the semi-direct
product of the other two groups).

\end{thm}

{\it Proof}. The proof will consist of two steps.

i) It is easy to see that the group $C_L^{\ast}$ of invertible elements in $C_L$ consists
of products $f(z)e_{\lambda}, \lambda \in L$ where $f(z)$ is an invertible analytic function.

 Let $W_L$ be a vector space consisting of $f\in C_L$ satisfying the following property:
for any invertible $a\in C_L$ there exists a constant $c$ such that  $afa^{-1} = f + c$.
Let us prove that such $f=f(z)$ is a linear function of $z$.
Indeed if $f=\sum_jf_j(z)e_{\lambda_i}$, then taking $a=e_{\mu}, \mu\in L$,
we see that the derivatives $f_j^{\prime}(z)$ are double-periodic analytic function on $\C$
(we do not assume that the periods are linearly independent over $\R$).
Hence $f_j(z)$ are linear functions in $z$.

Clearly $i(W_L)=W_{L^{\prime}}$.
 By composing $i$ with an
automorphism of $ C_L$ we can bring it to the
form $i(z)= \alpha z$ and $i(L)\subset L^{\prime}$, where $\alpha\in \C^{\ast}$.
Since $i$ is continuous we have $i(f(z))=f(\alpha z)$ for any $f\in {\cal O}(\C)$.

Notice that $e_{\lambda}ze_{-\lambda}=\lambda+z$. Since $e_{\lambda}$ is
invertible we have $i(e_{\lambda})=e_{\lambda^{\prime}}g(z)$ for some invertible
analytic function $g(z)$ and $\lambda^{\prime}\in L^{\prime}$.
Then $i(e_{\lambda}ze_{-\lambda})=\lambda+\alpha z$. On the other hand
it is equal to $\alpha \lambda^{\prime}+\alpha z$. It follows
that $\alpha L=L^{\prime}$.

ii) We let the groups $G_L$ and  $\C$ act on  $ C_L$ preserving the lattice $L$ and
sending $z \mapsto \alpha z$ , $\alpha \in G_L$ and $z \mapsto z+c$ ,
with $c\in \C$. It is easy to see that
the second automorphism is inner provided  $ c\in L$. It gives rise to an injective
 homomorphism  from
from the semi-direct product $G_L$ and $\C/L$ to  $Isom(E_L,E_{L})$.
  As we have just proven,
any automorphism of $C_L$ preserves the subspace of linear in $z$ functions.
This implies that the image of the latter homomorphism is a normal subgroup.
Hence, it is enough to identify the quotient with $Pic(E_L)$.

An element of the quotient
can be represented by an automorphism $i$ of $C_L$ identical on the space
of functions of $z$. We have $i(e_{\lambda})= f_{\lambda}e_{\lambda}$, where $f_{\lambda}$
is an invertible analytic function. It is easy to check that  $f_{\lambda}$ define a $1$-cocycle
of $L$ with coefficients in ${\cal O}^{\ast}$. If $i$ is an  inner automorphism,
the cocycle $f_{\lambda}$
is a coboundary. It gives rise to an  isomorphism  between the quotient and $H^1(L; {\cal O}^{\ast})=
Pic(E_L)$.  This proves the theorem. $\blacksquare$

\begin{cor} Two non-commutative elliptic curves $E_{\tau}$ and $E_{\mu}$ are
equivalent iff $\tau=g(\mu)$ for some $g\in  SL_2(\Z)$, where $SL_2(\Z)$ acts on
$\mu \in \C\setminus \Q$
by fractional transformations.
\end{cor}

It follows from the Corollary that the ``moduli space'' of non-commutative elliptic curves
can be identified with the quotient $\overline{D}/SL_2(\Z)$ where
$\overline{D}=\{z\in \C|\,|z|\le 1\}$ and $SL_2(\Z)$ acts by fractional transformations.
Naively, the moduli space consists of the isomorphism classes of algebras $C_L$ together
with the equivalence classes of the categories ${\cal B}_q$ where $q$ is a root of 1 (these
categories are equivalent to ${\cal B}_1$). The action of $SL_2(\Z)$ is not
discontinuous on the boundary circle $\partial D$. Hence the moduli space
is itself a non-commutative space. It can be described in terms of a cross-product algebra.
The reader should compare Corollary 3 with the similar theorem of Rieffel about
quantum tori (see [RS] for the higher-dimensional case).

\begin{rmk} It follows from the Corollary that if $\tau \in \Q(\sqrt{d})$,
where $d\in \Z$ is square-free then the
corresponding non-commutative elliptic curve has more autoequivalences than
a generic one.
If $d<0$ it is just an elliptic curve with complex multiplication. It was speculated
in [M2] that in the case $d>0$ non-commutative tori can be used for a description
of the maximal abelian extension of $\Q(\sqrt{d})$ (Stark conjectures).

\end{rmk}

Since our non-commutative elliptic curve can be defined not only over $\C$, the following
questions look very interesting.

\begin{que} Assume that we are working over the local field $K$ as in the Section $2$.

1. When two categories ${\cal B}_q$ and ${\cal B}_{q^{\prime}}$
are equivalent?

2. What is the group of auto-equivalences of ${\cal B}_q$?
\end{que}

We don't know the answer even in the case $K=\C, |q|=1$.

For usual elliptic curves (i.e. $|q|\ne 1$)  the answer is known thanks to the following result.
\begin{thm}
Assume that $|q|,|q^{\prime}| \ne 1$. Then

a) The categories ${\cal B}_q$ and ${\cal B}_{q^{\prime}}$
are equivalent if and only if the corresponding elliptic curves are
isomorphic.

b) The group of   auto-equivalences of ${\cal B}_q$ is generated
by the group of automorphisms of the elliptic
curve $E_q=\C^{\ast}/q^{\Z}$ and $Pic\, (E_q)$. The latter
 acts on  ${\cal B}_q$ sending a module $V$ to $V\otimes {\cal L}$,
where ${\cal L}\in Pic\, (E_q)$.
\end{thm}
{\it Proof.}\footnote{The proof bellow has been  explained to the authors by Dima Kaledin.}
Notice that irreducible objects of ${\cal B}_q$ are skyscraper sheaves supported at points of $E_q$.
Therefore, an equivalence $i: {\cal B}_q \to {\cal B}_{q^{\prime}}$ gives rise to a bijection $\tilde {i} $ between
the points of $E_q$ and $E_{q^{\prime}}$,
which is a homeomorphism with respect to the Zariski topology.
Further, line bundles ${\cal L}$ on $E_q$
can be characterised by the property that $dim(Hom({\cal L}, F))=1$,
for any skyscraper sheaf $F$. Hence  $i$ induces a bijection between line bundles on  $E_q$ and $E_{q^{\prime}}$.
 Therefore  it is enough to prove
that if $i({\cal O}_{E_q})={\cal O}_{E_q^{\prime} }$ then the functor $i$
is induced by an isomorphism of the ringed spaces.
Let $S\subset E_q(\C)$ be a finite subset. We are going to construct a morphism
of sheaves restricted to the complements of $S$
\begin{equation} \label{i}
i_{\ast}: {\cal O}(E_q \setminus S) \to {\cal O}(E_{q^{\prime}} \setminus \tilde {i}(S))
\end{equation}
Let $Coh_S(E_q)$ be the category of coherent sheaves  supported on $S$. The quotient category ${\cal B}_q/Coh_S(E_q)$
is equivalent to the category of coherent sheaves on $E_q\setminus S$. Furthermore  ${\cal O}(E_q \backslash S)$ can be identified
with the ring of endomorphism of the identity functor from $Coh(E_q\setminus S )$ to itself.

The functor $i$ induces an equivalence between  $Coh_S(E_q) $ and  $Coh_{ \tilde {i}(S)}(E_{q^{\prime}})$ which, in turn,
gives rise to an equivalence 
   $$i: Coh(E_q\setminus S ) \to Coh(E_{q^{\prime}}\setminus \tilde {i}( S) )$$
The latter defines the isomorphism (\ref{i}). $\blacksquare$

\begin{cor} If $L$ is a lattice then the group $Aut({\cal D}_L)$ of autoequivalences
of the category ${\cal D}_L$ is isomorphic to the group $Out(C_L)$ of outer automorphisms
of the algebra $C_L$.

\end{cor}

\begin{que} Is the Corollary $4$ true for a quasi-lattice $L$?

\end{que}

If the answer is positive, then we can say that two non-commutative elliptic curves
$E_L$ and $E_{L^{\prime}}$ are equivalent iff the correspodning categories
${\cal D}_L$ and ${\cal D}_{L^{\prime}}$ are equivalent. This observation agrees with the Theorem 6.

\subsection{Remark about higher genus curves}

If a complex curve $X$ has the genus $g\ge 2$ it admits a universal
covering by the upper-half plane $H$. We have $X\simeq H/\Gamma$,
where $\Gamma \subset SL_2(\Z)$ is a Fuchsian subgroup.
Applying Theorem 2 we see that $Coh(X)\simeq {\cal B}_{\Gamma}(H)$.
Now, if we consider a non-discrete embedding of $\Gamma$ into $SL_2(\Z)$,
we   define $Coh(H/\Gamma)$ as ${\cal B}_{\Gamma}(H)$.
Since the action of $\Gamma$ is now ergodic, we conclude that all objects of our category are free
${\cal O}(H)$-modules. In other words, there are no non-free coherent
sheaves on  degenerate curves. All the questions addressed above to  non-commutative
elliptic curves can be  addressed to non-commutative complex curves of
higher genus. Hopefully arising structures can be used for algebraic study of
Thurston boundary of Teichm\"uller space.

\section{Non-archimedean quantum tori}

This section is independent of the rest of the paper.

Let $K$ be a non-archimedean local field, $L$ is a free abelian
group of finite rank $d$, $\varphi:L\times L\to K$ is a skew-symmetric
bilinear form such that $|exp(2\pi i \varphi(\lambda,\mu))|=1$
for any $\lambda,\mu \in L$. We denote by
$A(T(L,\varphi))$ the
algebra of regular functions on the quantum torus $T(L,\varphi)$(see [So1]).
Thus it is an algebra with generators $e(\lambda),\lambda\in L$,
subject to the relation
$$e(\lambda)e(\mu)=exp(2\pi i \varphi(\lambda,\mu))e(\lambda+\mu).$$

There are various completions of  $A(T(L,\varphi))$ which appear
in the literature. In this section we will use the following one.

\begin{dfn} The space ${\cal O}(T(L,\varphi))$ of analytic functions
on the non-archime-\\
dean quantum torus consists of series
$\sum_{\lambda\in L}a(\lambda)e(\lambda),a(\lambda)\in K$ such that
$|a(\lambda)|\to 0$ as $|\lambda|\to \infty$ (here $|(\lambda_1,...,\lambda_d)|=
\sum_i|\lambda_i|$).

\end{dfn}

\begin{lmm} Analytic functions ${\cal O}(T(L,\varphi))$ form
a $K$-algebra.

\end{lmm}

{\it Proof.} We have $(\sum_{\lambda}a(\lambda)e(\lambda))
(\sum_{\mu}b(\mu)e(\mu))=\sum_{\nu}c(\nu)e(\nu)$,
where $c(\nu)=
\sum_{\lambda}a(\lambda)b(\nu-\lambda)exp(2\pi i \varphi(\lambda,\nu-\lambda))$.
We need to prove that \\
 $max_{\lambda,\nu\in L}|a(\lambda)||b(\nu-\lambda)|$
approaches zero as $|\nu|\to \infty$. Since both functions
$a(\lambda)$ and $b(\lambda)$ are bounded on $L$ we may assume that
$\lambda$ belongs to a finite subset of $L$. Then
$|a(\lambda)||b(\nu-\lambda)|\le const|b(\nu-\lambda)|$ where
$\lambda\in L_0$, and $L_0$ is a finite set. Therefore
$|c(\nu)|\to 0$ as $|\nu|\to \infty$. $\blacksquare$

Let us fix an isomorphism $L\simeq \Z^d$ and $r=(r_1,...,r_d)\in \R_{>0}^d$.
We denote by ${\cal O}(T(L,\varphi);r)$ the space of series
$f=\sum_{\lambda\in \Z^d}a(\lambda)e(\lambda)$ such that
$|a(\lambda)|r^{\lambda}\to 0$ as $|\lambda|\to \infty$.
Then ${\cal O}(T(L,\varphi);r)$ becomes a non-commutative Banach $K$-algebra
with the norm $||f||=max_{\lambda}|a(\lambda)|r^{\lambda}$.
In the case $\varphi=0$ an epimorphic image of ${\cal O}(T(L,\varphi);r)$
in the category of commutative Banach $K$-algebras is called
an affinoid $K$-algebra (see [Be]).
One can use non-archimedean quantum tori in order to define a deformation
of the sheaf of analytic functions on the rigid analytic space obtained from
degenerating hyper-K\"ahler manifolds in the spirit of  [KoSo1]. To be more precise,
such a degenerating family defines a fibration over an affine manifold $Y$ with
the fibers being tori (possibly degenerate on the codimension two subvariety
$Y^{sing}\subset Y$). The sheaf of analytic functions is in fact the sheaf
 on the non-singular part $Y\setminus Y^{sing}$ only. It is an open problem
to define it on $Y$.
One can replace series in commuting variables $x_i, 1\le i\le dim\,Y$ by the series
in variables $x_i$ such that $x_ix_j=qx_jx_i$. Then one obtains the deformation
of the sheaf ${\cal O}_Y$. Details will appear in the paper of M. Kontsevich
and the first author.

\vspace{15mm}

{\bf References}

\vspace{3mm}

[BG] V. Baranovsky, V. Ginzburg, Conjugacy classes in loop groups
and $G$-bundles on elliptic curves. Int. Math. Res. Notices, 15 (1996), 733-751.

\vspace{2mm}

[Be] V. Berkovich, Spectral theory and analytic geometry over \\
non-Archimedean fields. Mathematical Surveys and Monographs, 33.
Amer. Math. Soc., 1990.
\vspace{2mm}

[BDar] M. Bertolini, H. Darmon, The p-adic L-functions of modular
 elliptic curves, preprint, 2000.

\vspace{2mm}

\vspace{2mm}

[Bo] F. Boca, Projections in Rotation Algebras and Theta Functions,
math.OA/9803134

\vspace{2mm}

[BO] A. Bondal, D. Orlov,Reconstruction of a variety from
the derived category and groups of autoequivalences, alg-geom/97120029.

\vspace{2mm}

[CDS] A.Connes, M.Douglas, A. Schwarz, Non-commutative geometry
and Matrix theory: compactification on tori, hep-th/9711162.

\vspace{2mm}

[Co] A. Connes, Non-commutative geometry.

\vspace{2mm}

[CR] A. Connes, M. Rieffel, Yang-Mills for non-commutative two-tori.
Contemporary Math. 62 (1987), 237-266.

\vspace{2mm}

[De1] C. Deninger, On dynamical systems and their possible significance
for arithmetic geometry. In: Regulators in analyis, geometry and number
theory. Progress in mathematics 171, 29-87 (1999).

\vspace{2mm}

[De2] C. Deninger, Number theory and dynamical systems on foliated spaces,
preprint, 2001.

\vspace{2mm}

[Di] R. Dijkgraaf, Mirror symmetry and elliptic curves. In: Moduli
Space of Curves. Progress in Mathematics, vol. 129, p.149-163.

\vspace{2mm}

[Ga] P.  Gabriel,  Des categories abeliennes. Bull. Soc. Math. France 90 (1962)

\vspace{2mm}

 [GLO] V. Golyshev, V. Lunts, D. Orlov, Mirror symmetry for abelian varieties,
 math.AG/9812003.

\vspace{2mm}

[Ki] R. Kiehl, Theorem A und Theorem B in der nichtarchimedischen Funktionentheorie,
Inv. Math. 2 (1967), 256-273.

\vspace{2mm}

[Ko1] M. Kontsevich, Homological algebra of Mirror symmetry. Proceedings
of the ICM in Zurich, vol. 1, p. 120-139.

\vspace{2mm}

[Ko2] M. Kontsevich, Operads and Motives in Deformation Quantization,
 math.QA/9904055.

 \vspace{2mm}

[KR] M. Kontsevich, A. Rosenberg, Non-commutative smooth spaces,
math.AG/9812158.

\vspace{2mm}

[KoSo1]  M. Kontsevich, Y. Soibelman,
Homological mirror symmetry and torus fibrations,
math.SG/0011041, Published in Proceedings of KIAS conference
in Symplectic geometry and Strings, 2001.

\vspace{2mm}

[KoSo2]  M. Kontsevich, Y. Soibelman,  Deformation theory (book in preparation).

\vspace{2mm}

[M1] Yu. Manin, Quantized theta-function. Preprint RIMS, RIMS-700, 1990.

\vspace{2mm}

[M2] Yu. Manin,
Real multiplication and non-commutative geometry math.AG/0202109.

\vspace{2mm}

[M3] Yu. Manin, Mirror symmetry and quantization of abelian varieties,
math.AG/0005143.

\vspace{2mm}

[M4] Yu. Manin, The theory of commutative formal groups over fields of finite characteristic.
Russian Mathematical Surveys, Volume 18, n. 6.

\vspace{2mm}

[Mo1] G. Moore, Finite in all directions, hep-th/9305139.

\vspace{2mm}

[Mo2] G. Moore, Arithmetic and Attractors, hep-th/9807087

\vspace{2mm}

 [O1] D. Orlov,    On equivalences of derived categories of
 coherent sheaves on abelian varieties, alg-geom/9712017

\vspace{2mm}

[O2] D. Orlov, Equivalences of derived categories and K3 surfaces,
alg-geom/9606006.

\vspace{2mm}

[PZ] A. Polishchuk, E. Zaslow, Categorical mirror symmetry:
the elliptic curve, math.AG/9801119.

\vspace{2mm}

 [RS] M. Rieffel, A. Schwarz,
 Morita equivalence of multidimensional noncommutative tori, math.QA/9803057.

 \vspace{2mm}

[S1] A. Schwarz, Morita equivalence and duality, hep-th/9805034.

\vspace{2mm}

[S2] A. Schwarz, Theta-functions on noncommutative tori, QA/0107186.

\vspace{2mm}

[Sab] C. Sabbah, Systemes holonomes d'equations aux q-differences. In: D-modules
and microlocal geometry, 1992, p. 125-147.

\vspace{2mm}

[So1] Y. Soibelman, Quantum tori, mirror symmetry and deformation theory,
 math.QA/0011162

\vspace{2mm}

\vspace{3mm}

Addresses:

Y.S.:Department of Mathematics, KSU, Manhattan, KS 66506, USA

{email: soibel@math.ksu.edu}

V.V.: Department of Mathematics, University of Chicago, Chicago, IL 60637

{email: volgdsky@math.uchicago.edu}

\end{document}